\documentclass[a4paper]{article}
\usepackage[utf8]{inputenc}
\usepackage{amssymb,amsmath,latexsym}
\usepackage[english]{babel}

\usepackage{tikz-cd}

\title{On Field Generated by Division Points of Several Formal Groups}
\author{Soumyadip Sahu}
\date{8.8.2018}

\begin{document}
\maketitle
\begin{abstract}
In this article we prove some interesting results on field generated by division points of several formal groups, already implicit in the treatment in appendix-A of \cite{the}. More precisely, under suitable hypothesis we show that the field generated by the division points of several formal groups of same height is equal to the field generated by division points of an individual formal group among them, after a fixed finite unramified extension of base field.  
\end{abstract} 

\section{Introduction}
Let $p$ be an odd prime and let $K$ be a finite extension of $\mathbb{Q}_p$. Put $O_K$ to be the ring of integers of $K$, let $\mathfrak{p}_K$  
denote the unique maximal ideal of $O_K$, $v_{\mathfrak{p}_K}(\cdot)$ be the valuation associated to it and assume that the degree of the corresponding residue extension is $f$.\\ 
Fix an algebraic closure $\overline{\mathbb{Q}}_p$ and $|.|_p$ be a fixed extension of the absolute value. Let $\overline{O}$ be the ring of integers of $\overline{\mathbb{Q}}_p$ and $\overline{\mathfrak{p}}$ be the unique maximal ideal of $\overline{O}$. Clearly $\overline{\mathfrak{p}} \cap K = \mathfrak{p}_K$.\\
Fix a generator $\pi$ of $\mathfrak{p}_K$ and let $K_{R} = \mathbb{Q}_p(\pi)$. Use $R$ to denote the ring of integers of $K_R$, $\mathfrak{p}_R$ to denote the maximal ideal and $f_R$ to denote the degree of the residue extension.\\ 
Let $\mathfrak{F}$ be a (one dimensional, commutative) formal group law defined over $O_K$ which satisfies following additional conditions :\\~\\
i) $\mathfrak{F}$ has a formal $R$-module structure ,\\
ii) if \[[\pi](X) = \pi X + a_2X^2 +a_3X^3 +a_4X^4 +\cdots \in O_K[[X]] \] then $\min \, \{ i \geq 2 | |a_i|_{p} = 1 \} = p^h$ for some positive integer $h$. 
The integer $h$ will be called the height of group law. This condition is satisfied unless all the $a_i$-s are in maximal ideal (see \cite{haz2}, 18.3.1).\\~\\
A formal group law satisfying condition $(i)$ will be called a $\pi$-unramified group law.\\
Now $\mathfrak{F}$ defines a $R$-module structure on $\mathfrak{p}_K$ which naturally extends to a $R$-module structure on $\overline{\mathfrak{p}}$. We shall denote the corresponding addition by $\oplus_{\mathfrak{F}}$ to distinguish it from usual addition.\\
For each $n \geq 1$, let $\mathfrak{F}[\pi^n]$ denote the $\pi^n$-torsion submodule of $\overline{\mathfrak{p}}$ and let $L(\pi^n)$ be the subfield of $\overline{\mathbb{Q}}_p$ generated by $\mathfrak{F}[\pi^n]$ over $L$, where $L$ is any sub-field of $\overline{\mathbb{Q}_p}$. We shall adopt the convention $\mathfrak{F}[\pi^0] = \{0\}$.\\
Note that from the condition $(i)$ on $\mathfrak{F}$ implies \[[\pi^n](X) = \pi^n X + \text{higher degree terms} \tag{1.1}\] for each $n \geq 1$.\\
Note that if $K/\mathbb{Q}_p$ is an unramified extension then one can take $\pi = p$. In this case the condition $(i)$ on $\mathfrak{F}$ is verified trivially and $R = \mathbb{Z}_p$.\\
First, let us recall lemma-3.2 from appendix-A of \cite{the} :\\~\\
Let $\mathfrak{F}$ be a $\pi$-unramified group law of height $h$ defined over $O_K$ for some generator $\pi$ of $\mathfrak{p}_K$. Assume that $h | f$. Then $K(\pi)/K$ is a totally ramified cyclic extension of degree $(q-1)$.\\~\\
Now let $\pi_1$, $\pi_2$ be two generators of $\mathfrak{p}_K$, and let $\mathfrak{F}_1$ and $\mathfrak{F}_2$ be $\pi_1$-unramified and $\pi_2$-unramified formal group laws (resp.) of height $h_1, h_2$ (resp.) defined over $O_K$. The goal of this article is to describe the extension $K(\mathfrak{F}_1[\pi_1], \mathfrak{F}_2[\pi_2]) / K$ under the hypothesis $h_i | f$ for $i = 1,2$. At first we shall restrict to the case $h_1 = h_2$ and make some remarks about the general case later. \\~\\
The main result can be summerized as follows :\\~\\
\textbf{Theorem 1.1 :} Let $\mathfrak{F}_1, \cdots, \mathfrak{F}_n$ are $n$ unramified group laws of same height $h$ defined over $O_K$ corresponding to the generators (of $\mathfrak{p}_K$) $\pi_1, \cdots, \pi_n$ respectively, for some $n \geq 2$. Assume that $h| f$. Say, $K_1 = K(\mu_{(q-1)(p^f - 1)})$ where $q = p^h - 1$ . Then : \[K_1(\mathfrak{F}_1[\pi_1]) = \cdots = K_1(\mathfrak{F}_n[\pi_n]) = K_1(\mathfrak{F}_1[\pi_1], \cdots,\mathfrak{F}_n[\pi_n]).\]\\
We shall apply this theorem to deduce some results about the field generated by torsion points super-singular abelian varieties defined over an unramied extension of $\mathbb{Q}_p$. \\~\\
\textbf{Notations and conventions :}\\
 We shall use the notation 
\[
 q_h = p^h.
 \] We often use notations introduced in this section without explicitly defining them again.\\
A $p$-adic field is a finite extension of $\mathbb{Q}_p$.\\
For $n \geq 1$, $\mu_{n}$ be the group of $n$-th roots of unity in $\overline{\mathbb{Q}}_p$.\\
If $L$ is a $p$-adic field then $O_L$ is the ring of integers and $\mathfrak{p}_L$ is the unique maximal ideal of $O_L$ and residue degree of $L$ is $f_L$. Then $\mu_{p^{f_L} - 1} \subset O_L$ and $\mu_{p^{f_L} - 1} \cup \{0\}$ forms a cannonical set of representatives for the residue field. We shall always work with this set of representatives. \\
We shall often refer to appendix-A of \cite{the} just by mentioning appendix-A.

\section{Main Results}
Let $K$ be as in introduction. Let $\pi_1$, $\pi_2$ be two generators of $\mathfrak{p}_K$, the maximal ideal in ring of integers of $K$ and let $\mathfrak{F}_1$ and $\mathfrak{F}_2$ be $\pi_1$-unramified and $\pi_2$-unramified formal group laws respectively, of height $h_1$ and $h_2$ defined over $O_K$. The goal of this section is to describe the extension $K(\mathfrak{F}_1[\pi_1], \mathfrak{F}_2[\pi_2]) / K$ under the hypothesis $h_i | f$ for $i = 1, 2$. For simplicity we shall write $q_{h_i}$ as $q_{i}$  for all $i$ in range. Note that the hypothesis in preceding line means $\mu_{q_i - 1} \subseteq \mu_{p^f - 1} \subseteq O_K$ for $i = 1, 2$. \\~\\
From the proof of lemma-3.2 in appendix-A we conclude that there are pairs $(\Pi_i, \pi_i')$ such that $\Pi_i$ is a generator of $\mathfrak{p}_{K(\mathfrak{F}_i[\pi_i])}$, $\pi_i'$ is a generator of $\mathfrak{p}_{K}$, $\Pi_i^{q_i-1} = \pi_i'$ and $K(\Pi_i) = K(\mathfrak{F}_i[\pi_i])$ for $i = 1, 2$.\\
Fix such a pair for each $i = 1, 2$.\\
Consider the quotient $\frac{\pi_2'} {\pi_1'}$. Since degree of residue extension of $K$ is $f$ there is a $\zeta \in \mu_{p^f - 1}$ such that \[|\frac {\pi_2'} {\pi_1'} - \zeta|_p < 1 .\tag{2.1}\] Thus $\pi_2' = \zeta \pi_1' + x\pi_1'$ for some $x \in \mathfrak{p}_K$.\\
Let $\zeta_1$ be one $(q_1 - 1)$-th root of $\zeta$. With this notation we have the following lemma :\\~\\
\textbf{Lemma 2.1 :} Assume that $h_2 | h_1$. Then $ \Pi_2 \in K (\zeta_1\Pi_1)$.\\~\\
\textbf{Proof :} Proof of this lemma is similar to the proof of the lemma after proposition-12, chapter-2 in \cite{lang}. For sake of completion we repeat the argument.\\
Note that $\Pi_2^{q_2 - 1} = \pi_2'$. So it is enough to show that $X^{q_2-1} - \pi_2'$ has a root in $K(\zeta_1\Pi_1)$ since $\mu_{q_2-1} \subseteq K$. \\
Put $\beta = \zeta_1\Pi_1$. Then $\beta^{q_1-1} = \pi_2' + \pi_2'y$ where $y = - x \frac{\pi_1'} {\pi_2'} \in \mathfrak{p}_K$.\\
Let $f(X) = X^{q_1-1} - \pi_2'$ and let $\alpha_1, \cdots, \alpha_{q-1}$ be roots of $f$ in $\overline{\mathbb{Q}}_p$. Clearly \[|\alpha_1|_p = \cdots = |\alpha_{q_1-1}|_p = |\beta|_p = |\pi_2'|_{p}^{\frac {1} {q_1-1}}\]
Now \[f(\beta) = (\beta - \alpha_1)\cdots(\beta - \alpha_{q_1-1}).\]
But $f(\beta) = \pi_2'y$. So $|f(\beta)|_p < |\pi_2'|_p$. Hence there is at least one $i$ with $1 \leq i \leq q_1-1$ such that \[|\beta - \alpha_i|_p < |\pi_2'|_p^{\frac {1} {q_1-1}} = |\alpha_i|_p.\]
Fix such a $i$.
Note that \[|f'(\alpha_i)|_p = |\alpha_i|_p^{q_1-2} = \prod_{1 \leq j \leq q_1-1, j \neq i} |\alpha_j - \alpha_i|_p.\]
Since $|\alpha_j - \alpha_i|_p \leq |\alpha_i|_p$ for all $j \in \{1, \cdots, q_1-1\} - \{i\}$, we deduce that $|\alpha_j - \alpha_i|_p = |\alpha_i|_p$ for all $j$ in the range. This implies
\[  |\beta - \alpha_i|_p < |\alpha_j - \alpha_i|_p\] for all $j \in \{1, \cdots, q_1-1\} -\{i\}$.\\
Now Krasner's lemma (proposition-3, chapter-2, \cite{lang}) implies $\alpha_i \in K(\beta)$. \\
By assumption, $h_2 | h_1$. Hence $q_2 - 1 \, | \, q_1 - 1$. Put \[ \alpha = \alpha_i ^{\frac {q_1 - 1} {q_2 - 1}} \in K(\beta). \] Then $\alpha$ is a root of $X^{q_2 - 1} - \pi_2'.$ As argued before, this shows $\Pi_2 \in K(\zeta_1\Pi_1)$ as desired. $\square$ \\~\\
We shall treat $h_1 \neq h_2$ case in the following sub-section. For rest of this part we shall assume $h_2 = h_1 = h$. For simplicity put $q_1 = q_2 = q$.\\~\\
\textbf{Corollary 2.2 :} i) $K(\Pi_1, \Pi_2) = K(\zeta_1, \Pi_1)$.\\
ii) $K(\mathfrak{F}_1[\pi_1], \mathfrak{F}_2[\pi_2]) = K(\Pi_1, \Pi_2) = K(\zeta_1, \Pi_1) = K(\zeta_1, \Pi_2).$\\
iii) $K(\zeta_1)/K$ is the maximal unramified sub-extension of $K(\mathfrak{F}_1[\pi_1], \mathfrak{F}_2[\pi_2])/K$. \\~\\
\textbf{Proof :} $(i)$ By lemma 2.1, $\Pi_2 \in K(\zeta_1\Pi_1)$. Hence $K(\Pi_1, \Pi_2) \subseteq K(\zeta_1, \Pi_1)$. \\
Assume that the degree of residue extension corresponding to $K(\Pi_1, \Pi_2)$ is $f_1$. Thus there is a $\zeta_1' \in \mu_{p^{f_1} - 1} \subseteq K(\Pi_1, \Pi_2)$ satisfying $|\frac {\Pi_2} {\Pi_1} - \zeta_1'|_p < 1$. Ultrametric triangle inequality implies \[|\frac {\pi_2'} {\pi_1'} - \zeta_1'^{q-1}|_p < 1.\]
Comparing with $(2.1)$ and using ultrametric trinagle inequality we conclude that $|\zeta_1'^{q-1} - \zeta|_p < 1$. But $\zeta_1'^{q-1}, \zeta \in \mu_{p^{f_1} - 1}$. Thus $|\zeta_1'^{q-1} - \zeta|_p < 1$ implies $\zeta_1'^{q-1} = \zeta$. By definition of $\zeta_1$, $\zeta_1^{q-1}  = \zeta_1'^{q-1}$. Since $\mu_{q-1} \subseteq K$ and $\zeta_1' \in K(\Pi_1, \Pi_2)$, we have $\zeta_1 \in K(\Pi_1, \Pi_2)$. Hence $K(\zeta_1, \Pi_1) \subseteq K(\Pi_1, \Pi_2)$ and we have the result. $\square$\\~\\
$(ii)$ First equality follows from the facts that $K(\mathfrak{F}_1[\pi_1]) = K(\Pi_1)$ and $K(\mathfrak{F}_2[\pi_2]) = K(\Pi_2)$.\\  
Second equality is already proved in part $(i)$.\\
Imitation of proofs of lemma-2.1 and corollary-2.2 proves $K(\Pi_1, \Pi_2) = K(\zeta_1^{-1}, \Pi_2)$. Now the third equality follows from the fact $K(\zeta_1^{-1}, \Pi_2) = K(\zeta_1, \Pi_2)$. $\square$ \\~\\
$(iii)$ Note that $\zeta \in \mu_{p^f - 1}$ and $\zeta_1^{q-1} = \zeta$. Hence $\zeta_1 \in \mu_{(q-1)(p^f - 1)}$. Since $\gcd ((q-1)(p^f - 1), p) = 1$, the extension $K(\zeta_1)/K$ is unramified.\\
Now the extension $K(\Pi_1)/K$ is totally ramified (lemma-3.2, appendix-A).\\
So $K(\zeta_1) \cap K(\Pi_1) = K$ and $K(\zeta_1, \Pi_1)/K(\zeta_1)$ is a totally ramified extension (see lemma 1.2.1 in \cite{the} or lemma 2.1 in \cite{hab}). Thus $K(\zeta_1)/K$ is the maximal unramified sub-extension of $K(\zeta_1, \Pi_1)/K$. Now the result follows from the equality $K(\mathfrak{F}_1[\pi_1], \mathfrak{F}_2[\pi_2]) = K(\zeta_1, \Pi_1)$. $\square$\\~\\ 
\textbf{Corollary 2.3 :} Let $\zeta_1$ be as above. Then \[ \text{Gal}(K(\mathfrak{F}_{1}[\pi_1], \mathfrak{F}_{2}[\pi_2])/ K) \cong \text{Gal}(K(\zeta_1)|K) \times \mathbb{Z}/(q-1)\mathbb{Z}. \] \\
\textbf{Proof :} From lemma 1.2.1 of \cite{the} and the discussion in proof of corollary-2.2 it follows that \[ \text{Gal}(K(\mathfrak{F}_{1}[\pi_1], \mathfrak{F}_{2}[\pi_2])/ K) \cong \text{Gal}(K(\zeta_1)|K) \times \text{Gal}(K(\Pi_1)|K).\]
Now the result follows from \[\text{Gal}(K(\Pi_1)|K) \cong \mathbb{Z}/(q-1)\mathbb{Z}.\]
(See lemma-3.2 in appendix-A)\\
$\square$\\~\\
\textbf{Remark 2.4 :} Let $e$ be the smallest positive integer such that $\zeta_1^e \in K$. Clearly $e|(q-1)$ and $\mu_e \subseteq \mu_{q-1} \subseteq K$. Consider $f(X) = X^e - \zeta_1^e \in O_K[X]$. From definition of $f$, $f(\zeta_1) = 0$ and $K(\zeta_1)/K$ is an abelian Kummer extension of exponent $e$.\\~\\
\textbf{Corollary 2.5 :} Recall that $K_1 = K(\mu_{(q-1)(p^f-1)})$. With this notation we have :\\ \[K_1(\mathfrak{F}_1[\pi_1]) = K_1(\mathfrak{F}_2[\pi_2]).\]\\
\textbf{Proof :} Note that $K(\zeta_1, \Pi_1) \subseteq  K_1(\mathfrak{F}_1[\pi_1])$ and $K(\zeta_1, \Pi_2) \subseteq K_1(\mathfrak{F}_2[\pi_2])$. Now the result follows from corollary-2.2 $(ii)$. $\square$\\~\\
\textbf{Proof of theorem 1.1 :} Let $h$ be a positive integer and put $q = p^h$. Assume that $\mathfrak{F}_1, \cdots, \mathfrak{F}_n$ are $n$ unramified group laws of same height $h$ defined over $O_K$ corresponding to the generators (of $\mathfrak{p}_K$) $\pi_1, \cdots, \pi_n$ respectively, for some $n \geq 2$. We want to show : \[K_1(\mathfrak{F}_1[\pi_1]) = \cdots = K_1(\mathfrak{F}_n[\pi_n]) = K_1(\mathfrak{F}_1[\pi_1], \cdots ,\mathfrak{F}_n[\pi_n]).\]
It follows directly from corollary 2.5. $\square$\\~\\
\textbf{Remark 2.6 :} i) \[\text{Gal}(K_1(\mathfrak{F}_1[\pi_1], \cdots \mathfrak{F}_n[\pi_n])\,|\,K_1) = \text{Gal}(K_1(\mathfrak{F}_1[\pi_1]) | K_1) \cong \mathbb{Z}/(q-1)\mathbb{Z}.\]
This is a consequence of the facts that $K_1/K$ is unramified, $K(\mathfrak{F}_1[\pi_1])/K$ is totally ramified and $\text{Gal}(K(\mathfrak{F}_1[\pi_1])|K) \cong \mathbb{Z}/(q-1)\mathbb{Z} $.\\ 
ii) Consider the field $F= K(\mathfrak{F}_1[\pi_1], \cdots, \mathfrak{F}_{n}[\pi_n])$. From theorem-1.1 it follows that $F \subseteq K(\mu_{(q-1)(p^f - 1)}, \mathfrak{F}_i[\pi_i])$ for any $1 \leq i \leq n$. A more precise description of $F$ can be obtained as follows :\\
By corollary-2.2 for $1 \leq i < j \leq n$ there is a $\zeta_{ij} \in \mu_{(q-1)(p^f - 1)}$ such that $K(\mathfrak{F}_i[\pi_i], \mathfrak{F}_j[\pi_j]) = K(\zeta_{ij}, \Pi_i).$ Put $Z = \{ \zeta_{ij} \,|\, 1 \leq i < j \leq n \}.$ Then $F = K(Z, \mathfrak{F}_i[\pi_i])$ for any $1 \leq i \leq n$.\\
iii) Note that $K(Z)/K$ is unramified (since $Z \subseteq \mu_{(q-1)(p^f-1)}$) and $K(\mathfrak{F}_i[\pi_i])/K$ is tamely ramified (lemma-3.2, appendix-A). Hence $F/K$ is also tamely ramified.\\    
iv) One can give a precise description of $\zeta_{ij}$ in terms of $\mathfrak{F_i}$, $\mathfrak{F_j}$ and $\pi_i$, $\pi_j$.\\

\subsection{Unequal height case}
Assume that we are in the set-up as described in the beginning of the section. Further, $h_2 | h_1$. Let $(\Pi_1, \pi_1')$, $(\Pi_2, \pi_2')$, $\zeta$, $\zeta_1$ be as before. Clearly $K(\mathfrak{F}_1[\pi_1], \mathfrak{F}_2[\pi_2]) = K(\Pi_1, \Pi_2)$. Now we deduce a sequence of corollaries of lemma-2.1 as in previous part.\\~\\
\textbf{Corollary 2.1.1 :} $K(\Pi_1, \Pi_2) \subseteq K(\zeta_1, \Pi_1)$.\\~\\
\textbf{Proof :} Follows from lemma - 2.1. $\square$ \\~\\
Note that, here equality need not hold. Put $K_1 = K(\mu_{(q_1 - 1)(p^f - 1)})$. Then, \\~\\
\textbf{Corollary 2.1.2 :} $K_1(\mathfrak{F}_1[\pi_1], \mathfrak{F}_2[\pi_2]) = K_1(\mathfrak{F}_1[\pi_1])$.\\~\\
\textbf{Proof :} Follows from corollary-2.1.1 and the fact that $\zeta_1 \in \mu_{(q_1 - 1)(p^f - 1)}$. $\square$ \\~\\
One can prove an analogue of theorem 1.1 :\\~\\
\textbf{Proposition 2.1.3 :} Let $\mathfrak{F}_1, \cdots, \mathfrak{F}_n$ are $n$ unramified group laws of height $h_1, \cdots, h_n$ (resp.) defined over $O_K$ corresponding to the generators (of $\mathfrak{p}_K$) $\pi_1, \cdots, \pi_n$ (resp.), for some $n \geq 2$. Assume that $h_1| f$ and $h_i | h_1$ for all $2 \leq i \leq n$. Say, $K_1 = K(\mu_{(q_1-1)(p^f - 1)})$ where $q_1 = p^{h_1} - 1$ . Then : \[K_1(\mathfrak{F}_1[\pi_1]) = K_1(\mathfrak{F}_1[\pi_1], \cdots,\mathfrak{F}_n[\pi_n]).\]\\
\textbf{Proof :} Follows from corollary-2.1.2 . $\square$ \\~\\
\textbf{Remark 2.1.4 :} As before (with hypothesis in proposition-2.1.3) \[\text{Gal}(K_1(\mathfrak{F}_1[\pi_1], \cdots, \mathfrak{F}_n[\pi_n])|K_1) = \text{Gal}(K_1(\mathfrak{F}_1[\pi_1])|K_1) \cong \mathbb{Z}/(q_1-1)\mathbb{Z}.\]\\
\textbf{Corollary 2.1.5 :} With notation and hypothesis of the propositon above :\[K_1(\mathfrak{F}_1[\pi_2], \cdots,\mathfrak{F}_n[\pi_n]) \subseteq K_1(\mathfrak{F}_1[\pi_1]). \]\\
\textbf{Proof :} Follows from proposition-2.1.3. $\square$\\~\\
Corollary-2.1.5 gives an indication how one can describe the compositum if no divisibility relation among heights (as in hypothesis of proposition-2.1.3) is known. First, we need an existential result :\\~\\
\textbf{Proposition 2.1.6 :} Let $\pi$ be a generator of $\mathfrak{p}_K$ and use $A$ to denote the ring of integers of $\mathbb{Q}_p(\pi)$. Assume that the degree of the residue extension of $A$ is $f_A$. Let $h$ be a positive integer such that $f_A \,|\,h$. Then there is formal $A$-module of height $h$ defined over $O_K$.\\~\\
\textbf{Proof :} Follows from Hazewinkel's construction of universal formal $A$-module as in \cite{haz2}, 21.4.8. In particular, if one wants a formal $A$-module of height $nf_A$ for some $n \in \mathbb{N}$, then one have to (with notation in \cite{haz2}) map $S_2, \cdots, S_{nf_A - 1}$ to elements in $\mathfrak{p}_K$ and $S_{nf_A}$ to an unit in $O_K$.   $\square$\\~\\
Now let $\mathfrak{F}_1, \cdots, \mathfrak{F}_n$  be $n$ unramified group laws of height $h_1, \cdots, h_n$ (resp.) defined over $O_K$ corresponding to the generators $\pi_1, \cdots, \pi_n$ (resp.). Further, $h_i | f$ for all $1 \leq i \leq n$. Fix a generator $\pi$ (one may choose $\pi$ among $\pi_1, \cdots \pi_n$). Let $A$ and $f_A$ be as above. Put $H = \text{lcm} (h_1, \cdots, h_n, f_A)$ and assume that $\mathfrak{F}$ is a $\pi$-unramified group law defined over $O_K$. Note that $H\, |\, f$. Put $K_1 = K(\mu_{(p^H-1)(p^f - 1)})$. Now it follows from corollary-2.1.5 that \[K_1(\mathfrak{F}_[\pi_1], \cdots, \mathfrak{F}_n[\pi_n]) \subseteq K_1(\mathfrak{F}[\pi]).\]

\section {An application}
Let $\mathbb{Q}_{w}$ be an unramified extension of $\mathbb{Q}_p$ such that the cardinality of field of residues associated to $\mathbb{Q}_w$ is $w = p^f$ for some integer $f \geq 1$. Denote the ring of integers of $\mathbb{Q}_w$ to be $\mathbb{Z}_w$. Assume that $\mathfrak{F}_1, \cdots, \mathfrak{F}_n$ are $n$ formal group laws defined over $\mathbb{Z}_w$ of same height $h$ for some $n \geq 1$ and $F_1, \cdots, F_n$ are the formal power series corresponding to $\mathfrak{F}_1, \cdots, \mathfrak{F}_n$ respectively. Further, $h\,|\,f$.\\
Consider a tuple \[F(X_1, \cdots, X_n, Y_1 \cdots, Y_n) = (F_1(X_1, Y_1), \cdots, F_n(X_n, Y_n)).\]
It defines a n-dimentional formal group law over $\mathbb{Z}_w$ (say $\mathfrak{F}$). (See \cite{haz2}, II.9.1)\\
$\mathfrak{F}$ defines a $\mathbb{Z}_p$ module structure on $(\overline{\mathfrak{p}})^{\times n}$. Note that the addition and scalar multiplication is component-wise. \\
Let $\mathfrak{F}[p]$ denote the subgroup of $p$-torsion points of this module. For any subfield $L \subseteq \overline{\mathbb{Q}}_p$, let $L(\mathfrak{F}[p])$ be the field generated over $L$ by co-ordinates of $p$-torsion points of of $\mathfrak{F}$. Note that $L(\mathfrak{F}[p]) = L(\mathfrak{F}_1[p], \cdots, \mathfrak{F}_n[p])$.\\
Put $L = \mathbb{Q}_w(\mu_{(q_h - 1)(w-1)})$. Then :\\~\\
\textbf{Proposition 3.1 :} i) $L(\mathfrak{F}[p]) = L(\mathfrak{F}_i[p])$ for any $1 \leq i \leq n$.\\
ii) The extension $\mathbb{Q}_w(\mathfrak{F}[p])/\mathbb{Q}_w$ is tamely ramified.\\~\\
\textbf{Proof :} $(i)$ Follows from theorem-1.1.\\
$(ii)$ Follows from remark-2.6 $(ii)$. $\square$\\~\\
One can generalize proposition-3.1 as theorem-1.1 can be generalized to proposition-2.1.3, ie one can start with $n$ group laws of distict height satisfying $h_i | h_1 \, \forall\,1 \leq i \leq n$, $h_1 | f$ and write down an analogue of propositon-2.1.3. \\~\\
Now let $A$ be an abelian variety which is defined over $\mathbb{Q}_w$ and is isomorphic to a product of super-singular elliptic curves over $\mathbb{Q}_w$. Note that all super-singular elliptic curves define group law of height $2$ and all $p$-torsion points of a super-singular elliptic curve comes from the formal group associated to it. Assume $2 | f$. Let $A[p]$ denote the group of $p$-torsion points of $A$. Then :\\~\\
\textbf{Proposition 3.2 :} $\mathbb{Q}_w(A[p])$ is a tamely ramified extension of $\mathbb{Q}_w$.\\~\\
\textbf{Proof :} Follows from proposition - 3.1. $\square$\\~\\
\textbf{Remark 3.3 :} Let $A$ be any super-singular abelian variety defined over an unramified extension of $\mathbb{Q}_p$. This means $A$ has a good reduction at $p$ and the reduced abelian variety $\widetilde{A}$ is super-singular. Now it is well-known that $\widetilde{A}$ is isogenous to a product of super-singular elliptic curves in some finite extension of field of definition of $\widetilde{A}$ (Oort). If this isogeny can be lifted to some finite unramified extension of $\mathbb{Q}_p$ and if its degree is not divisible by $p$, we can possibly make a statement analogus to proposition-3.2 for arbitrary super-singular abelian varieties.

Soumyadip Sahu\\
Kolkata, India.\\
soumyadip.sahu00@gmail.com
\end{document}